%

%
%
\documentstyle{amsppt}
\define\stick{
\setbox255=\hbox{\raise0.6ex\hbox{$\scriptstyle\bullet$}}
\mathord{\rlap{\hbox to\wd255{\hss\hbox{$\scriptstyle|$}\hss}}
\box255}
}

\define\restr%
{{\hskip-0.1ex|\hskip-0.02ex{\grave{}}\hskip0.5ex}}
\def\dom{\hbox{dom}}

\topmatter

\title \nofrills  Luzin and anti-Luzin almost disjoint families
\endtitle
\author  Judith Roitman${}^1$
\thanks {${}^1$ Partially supported by
the University of  Kansas General Research Fund grant \# 3058}
\endthanks
and Lajos Soukup${}^2$
\thanks {${}^2$ Supported by the Hungarian National Foundation for Scientific
Research grant \# 16391   } \endthanks
\endauthor
\affil  University of Kansas, 
Mathematical Institute of the Hungarian Academy
of Science            \endaffil
\address  Lawrence, KS 66045; POB 127, H-1364 Budapest, Hungary
\endaddress 
\email roitman\@math.ukans.edu, soukup\@math-inst.hu
\endemail
\dedicatory
\enddedicatory
\date   November, 1997         \enddate
\abstract  Under $MA_{\omega_{1}}$ every uncountable almost disjoint family is
either anti-Luzin or has an uncountable Luzin subfamily.  This fails under
CH.  Related properties are also investigated.
\endabstract
\leftheadtext{Judith Roitman and Lajos Soukup}
\endtopmatter

\noindent {\bf 0. Introduction}
\medskip
This paper looks at two combinatorial
properties of almost disjoint families,
Luzin and anti-Luzin, along with some variants.

Let $\Cal A$ be an uncountable almost disjoint family
on a countable set $W$.

\definition {Definition 0.1}
$\Cal A$ is Luzin iff there is an enumeration $\Cal A =
\{a_{\alpha}:
\alpha <
\omega_{1}\}$ such that for all $w \in [W]^{<\omega}$, $\alpha < \omega_{1}$,
$\{\beta <
\alpha: a_{\alpha} \cap a_{\beta} \subset w\}$ is finite.
\enddefinition

\definition {Definition 0.2}  $\Cal A$ is anti-Luzin iff $\forall \Cal B \in
[\Cal A]^{\omega_{1}}
\exists
\Cal C,
\Cal D \in [\Cal B]^{\omega_{1}} \ \bigcup \Cal C \cap \bigcup \Cal D$ is
finite.
\enddefinition

The notion of ``Luzin" is somewhat standard.  The notion of ``anti-Luzin" is
new.

Luzin almost disjoint families are an analogue of Hausdorff gaps: they exist
(via a diagonal construction) in ZFC, and because of the finitary nature of the
definition the property of ``Luzin" is upwards absolute.  Perhaps because of
this similarity, they have been called Luzin {\it gaps}.

Anti-Luzin families also exist in ZFC, as canonical objects:  any uncountable
set of branches of a countable tree $T$ is an anti-Luzin family on $T$.

Luzin and anti-Luzin families are hereditary in the following sense:
uncountable
subfamilies of Luzin (respectively anti-Luzin) almost disjoint families are
Luzin (respectively anti-Luzin).

\definition {Definition 0.3}  $\Cal A$ is near-Luzin iff $\forall \Cal C, \Cal
D
\in [\Cal A]^{\omega_1} \ \bigcup \Cal C \cap \bigcup \Cal D$ is infinite.
\enddefinition

We will rely heavily on the following obvious statement: 
$\Cal A$ is anti-Luzin iff it has no uncountable near-Luzin subfamily.

Near-Luzin first appeared in \cite{HJ}, where it was called
$\omega_1$-full: near-Luzin families give rise to compact Hausdorff spaces in
which the intersection of any two uncountable open sets is non-empty.

\proclaim {Claim 0.4}  If $\Cal A$ is Luzin then it is near-Luzin.
\endproclaim

\demo {Proof}  Fix an arbitrary enumeration $\{a_{\alpha}: \alpha <
\omega_1\}$ of $\Cal A$.  Suppose $\Cal A$ is not near-Luzin.  Then there are
uncountable $\Cal C, \Cal D \subset \Cal A$ with $\bigcup \Cal C \cap \bigcup
\Cal D = w$ finite.  So for each $c \in \Cal C, d \in \Cal D, c \cap d
\subset w$.  There is $\alpha$ with $a_{\alpha} \in \Cal C, \{a_{\beta} \in
\Cal D: \beta < \alpha\}$ infinite.  So $\Cal A$ is not Luzin.
\enddemo

\proclaim {Corollary 0.5} An anti-Luzin family is not Luzin.
\endproclaim

While superficially corollary 0.5 does not reverse (a disjoint union of a Luzin
and an anti-Luzin family is neither Luzin nor anti-Luzin), does it reverse
in any
deep sense?  In particular, must any uncountable almost disjoint family which
does not embed one embed the other?  The answer is yes and no.

\proclaim {Theorem 0.6}  Assume MA$_{\omega_{1}}$.  Every uncountable almost
disjoint family is either anti-Luzin or contains an uncountable Luzin
subfamily.
\endproclaim

\proclaim {Theorem 0.7}  Assume $\stick$.  There is an uncountable almost
disjoint family which contains no uncountable anti-Luzin and no uncountable
Luzin subfamilies.
\endproclaim

Here $\stick$ is the following weakening of CH:  There is a family $\Cal S
\subset [\omega_1]^{\omega}$ of size
$\omega_1$ so that every uncountable subset of $\omega_1$ contains a set in
$\Cal S$.

Theorem 0.6 says that under MA$_{\omega_{1}}$, almost disjoint families have a
lot of structure.  Theorem 0.7 says that under $\stick$ they don't.
This used to be what one would expect, but recent work on iterating totally
proper forcing while preserving CH has changed our expectations.  In
particular, Abraham and Todorcevic showed the consistency of ``CH + all
$(\omega_{1},\omega_{1})$-gaps contain an uncountable Hausdorff sub-gap."
If Luzin almost disjoint families were combinatorially similar to
Hausdorff gaps the conclusion of theorem 0.6 would also be consistent with CH.
Thus theorem 0.7 destroys the parallel between Luzin almost disjoint families
and Hausdorff gaps.

In section 1 we prove theorem 0.6, in section 2 we prove theorem 0.7, in the
rest of the paper we explore some of the fine combinatorial structure of these
notions.

\definition {Conventions}  In this paper almost disjoint families are
collections of infinite sets whose pairwise intersections are finite; the
superscript ``$^*$" means ``mod finite; " all trees grow upward; and
properties are listed consecutively no matter what theorem, lemma, or
definition they occur in, so if there is a reference to property 17 the reader
can easily find it.
\enddefinition

\definition{Acknowledgment}
We would like to express our thanks to the referees for
 his (her) numerous valuable suggestions and comments.
\enddefinition
\bigskip

\noindent {\bf 1.  Proof of Theorem 0.6}
\medskip

Assume MA$_{\omega_{1}}$.  Let $\Cal A = \{a_{\alpha}: \alpha <
\omega_{1}\}$ be
an almost disjoint family which is not anti-Luzin.  By moving to a possibly
smaller subfamily, we may assume $\Cal A$ is near-Luzin.

Let ${\Bbb P}$ be the set $[\omega_{1}]^{<\omega}$ under the following partial
order:  $$p \leq q \ \roman{iff}\  p \supset q \ \roman{and} \ \forall
\alpha \in
p
\setminus q \forall \beta \in q \ \roman{if} \ \alpha < \beta \ \roman{then} \
a_{\alpha}
\cap a_{\beta} \not\subset k_{q}$$ where $k_{q} =$ max $\bigcup\{a_{\alpha}
\cap a_{\beta}: \{\alpha,\beta\} \in [q]^{2}\}$.
\medskip

\proclaim {Remark 1.1}  If $\beta >$ sup $p$ then $p \cup \{\beta\} \leq p$.
\endproclaim

A quick ad hoc definition: for $s, t \subset$ ON, $s \ll t$ iff sup $s <$ inf
$t$.

\proclaim {Remark 1.2}  If $r,
s, t$ are disjoint, $r \ll s \cup t$, $p = r \cup s$, $q = r
\cup t$ and $k_{p} = k_{q} = k$, then $p, q$ are compatible iff for all $\alpha
\in s, \beta \in t \ a_{\alpha} \cap a_{\beta} \not\subset k$.
\endproclaim

By remark 1.1, enough sets are dense so that if ${\Bbb P}$ is ccc
then, by MA$_{\omega_{1}}$, there is a generic filter $G$ so that
$\{a_{\alpha}:
\alpha
\in
\bigcup G\}$ is Luzin.  Remark 1.2 will be used to show ccc.
\medskip

\proclaim {Lemma 1.3}  Suppose $n, k \in \omega$ and suppose $E \in
\bigl[[\omega_{1}]^{n}\bigr]^{\omega_{1}}$ is pairwise disjoint.
Then there are $s \neq t \in E$ with
$a_{\alpha} \cap a_{\beta} \not\subset k$ for all $\alpha \in s, \beta \in t$.
\endproclaim

To prove lemma 1.3 we need

\proclaim {Sublemma 1.3.1}  If $S, T$ are uncountable disjoint subsets of
$\omega_{1}$ and
$k <
\omega$ then there are $S' \in [S]^{\omega_{1}}$  $T' \in [T]^{\omega_{1}}$ and
$k' > k$ with $k' \in a_{\alpha} \cap a_{\beta}$ for all $\alpha \in S'$,
$\beta \in T'$.
\endproclaim

\demo {Proof of sublemma 1.3.1}  Otherwise for all $k' > k$ either $S(k') =
\{\alpha \in S: k' \in a_{\alpha}\}$ is countable or $T(k') = \{\alpha \in
T: k'
\in a_{\alpha}\}$ is countable.  Let $S'' = S \setminus \bigcup \{S(k'): k'
> k$ and $S(k')$ is countable$\}$.  Let $T'' = T \setminus \bigcup \{T(k'):
k' > k$ and $T(k')$ is countable$\}$.  $S''$ and $T''$ are uncountable and
$\bigcup_{\alpha
\in S''}(a_{\alpha}
\setminus k) \cap \bigcup_{\beta \in T''}(a_{\beta}
\setminus k) = \emptyset$, contradicting our hypothesis on $\Cal A$.
\enddemo

\demo {Proof of lemma 1.3}  For $s \in E$ we write $s = \{\alpha(s,i): i <
n\}$  Let
$\{(i_{m}, j_{m}): m < n^{2}\}$ enumerate
$n \times n$.  Using lemma 1.3.1 iteratively, at each stage $m$ we have
uncountable disjoint subsets of $E$, $E(m)$,
$F(m)$, with $E(m) \supset E(m+1), F(m) \supset F(m+1)$, and
$k_{m} > k$ with $k_{m} \in a_{\alpha(s,i_{m})} \cap a_{\alpha(t,j_{m})}$ for
all $s \in E(m), t \in F(m)$.  But then for $\alpha \in s \in E(n^2-1), \beta
\in t \in F(n^2-1), a_{\alpha} \cap a_{\beta} \not\subset k$.
\enddemo

\proclaim {Lemma 1.4}  ${\Bbb P}$ is ccc.
\endproclaim

\demo {Proof}  Given an uncountable subset
$F$ of ${\Bbb P}$ we may without loss of generality assume that for some $n, F
\subset [\omega_{1}]^{n}$; $F$ is a
$\Delta$-system with root $r$ so that each $p \in F$ has the form $r \cup
s_{p}$; there is some $k$ with $\forall p \in F \ k = k_{p}$; and
$E =
\{s_{p}: p
\in F\}$ is well-ordered by
$\ll$.  By the lemma there are $p \neq q \in F$ so that for all
$\alpha \in s_{p}, \beta \in s_{q}$, $a_{\alpha} \cap a_{\beta} \not\subset
k$.
By remark 1.2, $p, q$ are compatible.
\enddemo
\bigskip

\noindent {\bf 2.  Proof of Theorem 0.7}
\medskip

In this section we prove theorem 0.7.  In the next section we will give a
stronger version, but the combinatorics are sufficiently complicated that it
makes sense to give the weaker proof first and then show how to improve it.

Let $\{S_{\alpha}: \alpha < \omega_1\}$ be a $\stick$-sequence, i.e.
each $S_{\alpha} \subset \alpha$ and every uncountable $X \subset \omega_1$
contains some
$S_{\alpha}$.

The
family
$\Cal A =
\{f_{\alpha}:
\alpha <
\omega_1\}$ will be a subset of
$\wp(\omega\times\omega)$, where each $f_{\alpha}$ is a function from
$\omega$ to $\omega$.

We require:

\hskip .3 in {\vbox{\hsize 4.5 in

\flushpar\hangindent=.8cm  1.  if $\alpha \neq \gamma$ then $\{i:
f_{\alpha}(i) =
f_{\gamma}(i)\}$ is finite

\flushpar\hangindent=.8cm  2.  if $\beta \leq \alpha$ then $f_{\alpha} \cap
\bigcup_{\gamma
\in S_{\beta}}f_{\gamma}$ is infinite

\flushpar\hangindent=.8cm  3.  if $\beta \leq \alpha$ then $\exists n \
\{\gamma
\in S_{\beta}: f_{\alpha} \cap f_{\gamma} \subset n\times n\}$ is infinite.

}}

Property 1 makes $\Cal A$ almost disjoint and property 2 makes it near-Luzin
(hence prevents an uncountable anti-Luzin subfamily). Finally, property 3
prevents an
uncountable Luzin subfamily. Indeed, let 
$\{f_{{\alpha}_{\nu}}:{\nu}<{\omega}_1\}$ be an enumeration of an uncountable
subfamily of $\Cal A$. Pick ${\beta}<{\omega}_1$ with
$S_{\beta}\subset \{{\alpha}_{\nu}:{\nu}<{\omega}_1\}$.
Let $I=\{{\nu}:{\alpha}_{\nu}\in S_{\beta}\}$ and fix
${\mu}<{\omega}_1$ such that ${\alpha}_{\mu}\ge {\beta}$ and
${\mu}>\sup I$. Then, by 3, there is $m\in {\omega}$ such that 
$
\{{\nu}<{\mu}:f_{{\alpha}_{\nu}}\cap f_{{\alpha}_{\mu}}\subset m\times m\}
$ is infinite, which contradicts the requirement formulated in
definition 0.1. So $\{f_{{\alpha}_{\nu}}:{\nu}<{\omega}_1\}$ is not Luzin.

Some preliminaries:

\definition {Definition 2.1} A set $F \subset \omega \times \omega$ is fat
iff $\limsup_n|\pi_n F| = \omega$, where $\pi_n F = \{j: (n,j) \in F\}$.
\enddefinition

Notice that if $\exists U$ a finite family of functions with $F \subset^*
\bigcup U$ then $F$ is not fat, and that fat sets are infinite.

\proclaim {Lemma 2.2} If $\Cal C$ is an infinite almost disjoint family of
functions from $\omega$ to $\omega$, then
$\bigcup
\Cal C$ is fat.
\endproclaim

\demo {Proof}  Fix $m <\omega$.  Let $G \in [\Cal C]^{m+1}$.  There is $r$
so that if $f \neq g \in G$ then $f{\restr}{[r,\omega)} \cap g{\restr}{[r,\omega)} =
\emptyset$.  So $m <$ lim sup$_n |\pi_n\bigcup G|$.
\enddemo

\definition {Definition 2.3}  Let $\Cal C$ be a countable almost disjoint
collection of functions from $\omega$ to $\omega$.  A finite partial function
$\sigma$ from $\omega$ to $\omega$ is $\Cal C$-free iff $\exists^{\infty}g \in
\Cal C \ g
\cap
\sigma =
\emptyset$.
\enddefinition

\proclaim {Lemma 2.4} Let $\Cal C$ be a countable collection of
functions from $\omega$ to $\omega$.  There is a function $s:
\omega \rightarrow \omega+1$ so that if $\sigma$ is a finite partial function
from $\omega$ to $\omega$ with $\sigma(i) \neq s(i)$ for all $i
\in$ dom
$\sigma$, then $\sigma$ is $\Cal C$-free.
\endproclaim

Such an $s$ is called $\Cal C$-tight.

\demo {Proof}  Consider $\Cal C$ as a subset of the compact space
$(\omega+1)^{\omega}$.  Let $s$ be an accumulation point of $\Cal C$ in
$(\omega+1)^{\omega}$.  If $\sigma$ is a finite partial function from
$\omega$ to
$\omega$ with $s \cap \sigma = \emptyset$ then $s \notin \{f:
f \cap
\sigma \neq \emptyset\}$ which is closed, so $\Cal C \setminus \{f:
f \cap
\sigma \neq \emptyset\}$ is
infinite.

\enddemo

The construction is a straightforward induction, given the following

\proclaim {Lemma 2.5}  Let $\Cal C$ be a countable almost disjoint collection
of functions from $\omega$ to $\omega$, let $\Cal C_n \subset \Cal C$ for
each $n < \omega$, and let $F_n$ be fat for each $n < \omega$.
Then there is a function
$f: \omega \rightarrow \omega$ so

\hskip .3 in {\vbox{\hsize 4.5 in

\flushpar\hangindent=.8cm  4.  $\{f\} \cup \Cal C$ is almost disjoint

\flushpar\hangindent=.8cm  5.  for each $n$, $f \cap
F_n$ is infinite

\flushpar\hangindent=.8cm  6.  $\forall n \exists m_n \ \{g \in \Cal C_n: f
\cap g \subset {m_n} \times \{f(i): i < m_n\}\}$ is infinite.

}}

\endproclaim


The family $\{f_{\alpha}:{\alpha}<{\omega}_1\}$ will be constructed
recursively in ${\omega}_1$ steps.
Assume that $\{f_{\alpha}:{\alpha}<{\beta}\}$ is already constructed.
Fix an enumeration $\{{\beta}_n:n<{\omega}\}$ of ${\beta}$. Let
$\Cal C=\{f_{\alpha}:{\alpha}<{\beta}\}$,
$\Cal C_n=\{f_{\alpha}:{\alpha}\in S_{{\beta}_n}\}$ and
$F_n=\cup\Cal C_n$. Now we can apply lemma 2.5 to get $f_{\beta}$ as $f$.
  
In
the next section, we will need to deal with many more fat sets, which
is why 2.5 is stated in its current generality.

\demo  {Proof of lemma 2.5}
Let $\Cal C = \{g_i: i < \omega\}$.

At stage $j$ we construct a
finite set $U_j$ (``$U$" is short for ``used up") of functions in $\Cal C$
where
$U_{j-1}
\subset U_j$; we define $m_j \geq j$, and define $f$ on
$(m_{j-1},m_j]$.

So suppose we are at stage $j$.  We know $m_k$ for each $k < j$,
$f{\restr}{m_{j-1}+1}, U_{j-1}$, and, for each $k < j$, we have an $\Cal
C_k$-tight
$s_k$.  Our induction hypothesis is that
$$\forall k < j, f{\restr}{(m_{k-1},m_{j-1}]} \cap s_k{\restr}{(m_{k-1},m_{j-1}]} =
\emptyset.$$

Let
$s_j$ be $\Cal C_j$-tight.

Since $F_k$ is fat, $\forall k \leq j$ there is some $r_{k,j} > m_{j-1}$ with
each
$r_{k,j} < r_{k+1,j}$, and some $(r_{k,j},t_{k,j}) \in F_k \setminus \bigcup
U_{j-1}$.  Let
$m_j = r_{j,j}$.  Let $r_{-1,j} = m_{j-1}$.

For $k \leq j$, let $f(r_{k,j}) = t_{k,j}$.  (This is towards
property 5.)

For
$i
\in (r_{-1,j},r_{j,j})$ with $\forall k \ i \neq r_{k,j}$, let $f(i)$ be any $m
\notin
\{g(i): i \in U_{j-1}\} \cup \{s_k(i): k \leq j\}$.  (This is towards
properties 4 and 6.)

For each $k \leq j$ $$f{\restr}{(m_{k-1},m_j]} \cap s_k{\restr}{(m_{k-1},m_j]} =
\emptyset$$ so
$f{\restr}{(m_{k-1},m_j]}$ is $\Cal C_k$-free, for each $k \leq j$.  Hence, for
each $k \leq j$ there is
$g_{k,j} \in \Cal C_k \setminus U_{j-1}$ with $$g_{k,j} \cap
f{\restr}{(m_{k-1},m_j]} =
\emptyset.$$

Let $$U_j = U_{j-1} \cup \{g_{k,j}: k \leq j\} \cup
\{g_j\}.$$

Property 4 is satisfied:  If $n \leq j$ and $i \in (m_{j-1},m_j]$ then
$g_n(i) \neq f(i)$, so $$f \cap g_n \subset m_n \times \{f(i): i < m_n\}.$$

Property 5 is satisfied:  $$\forall j \geq n \ f(r_{n,j}) = t_{n,j}$$
and $$(r_{n,j},t_{n,j}) \in F_n$$ so $f \cap F_n$
is infinite.

Property 6 is satisfied:  If
$j \geq n$ then $$f{\restr}{(m_{n-1},m_j]} \cap g_{n,j} = \emptyset$$ so $$g_{n,j}
\cap f
\subset m_n \times \{f(i): i < m_n\}.$$
\enddemo
\bigskip

\noindent {\bf 3.  A strengthening of theorem 0.7}
\medskip

In this section we strengthen theorem 0.7.

\definition {Definition 3.1}  An uncountable almost disjoint family $\Cal A$ is
strongly near-Luzin iff,  for every $\Cal C_0, ... \Cal C_n \in
[\Cal A]^{\omega_1} , \ \bigcap_{i \leq n}\bigcup \Cal C_i$ is infinite.
\enddefinition

Strongly near-Luzin families appear in \cite{JN}, where they are called strong Luzin families.  They cannot
exist under MA + $\neg$ CH.  The following theorem shows that they need not be
Luzin.

\proclaim {Theorem 3.2}  Assume $\stick$.  There is an uncountable almost
disjoint family which is strongly near-Luzin, but has no uncountable Luzin
subfamilies.
\endproclaim

The proof is somewhat like that of theorem 0.7, but the combinatorics
are more complicated, so complicated that we will invoke
elementary submodels to avoid stating them explicitly.

So let $\{S_{\alpha}: \alpha < \omega_1\}$ be a $\stick$-sequence.

We begin by strengthening property 2 to

\indent\indent 7.  If $\beta_0, ... \beta_n < \alpha$ and each sup
$S_{\beta_i} \cup (\beta_i+1) <$  $\inf S_{\beta_{i+1}}$ then $f_\alpha \cap
\bigcap_{i \leq n}
\bigcup_{\gamma \in S_{\beta_i}}f_{\gamma}$ is infinite.

We will be done if the sequence of $f_{\alpha}$'s satisfies properties 1,
7, and 3. Indeed, as we have seen in the proof of theorem 0.7,
property 3 implies that $\Cal A$  does not contain an uncountable
Luzin subfamily. Property 1 yields that $\Cal A$ is almost disjoint.
So we need to show that if property 7 holds then $\Cal A$ is strongly
near Luzin. So let $n\in {\omega}$ and 
$I_0,\dots, I_{n-1}\in [{\omega}_1]^{{\omega}_1}$. Since 
$\{S_{\alpha}:{\alpha}<{\omega}_1\}$ is a $\stick$-sequence we can
find ${\beta}_0<{\beta}_1<\dots<{\beta}_{n-1}\in {\omega}_1$ such that
$S_{{\beta}_i}\subset I_i$ and 
$\sup S_{{\beta}_{i-1}}\cup ({\beta}_{i-1}+1)<\min S_{{\beta}_i}$.
Then, by property 7, 
$\bigcap_{i<n}\bigcup\{f_{\gamma}:{\gamma}\in I_i\}\supset
\bigcap_{i<n}\bigcup\{f_{\gamma}:{\gamma}\in S_{{\beta}_i}\}$ is
infinite,
which was to be proved.

In applying 2.5 in the previous section we had the luxury of knowing that
each $\bigcup \Cal C_n$ was fat.  But an intersection of fat sets need not
be fat.  So we must ensure that the following property holds:

\indent\indent 8.  If $\beta_0, ... \beta_n <
\alpha$ and each sup
$S_{\beta_i} \cup (\beta_i+1) <$  $\inf S_{\beta_{i+1}}$, then $\bigcap_{i \leq
n}
\bigcup_{\gamma \in S_{\beta_i}}f_{\gamma}$ is fat.

Property 8 allows us to construct a family in which property 7 holds.  How will
we build a family in which property 8 holds?

\definition {Definition 3.4}  For $F \subset \omega \times \omega$ and
$E \in[\omega]^{\omega}$ let   $\pi_EF = F \cap (E \times \omega)$.
\enddefinition

To get property 8 to hold, we need to start with enough fat $C$'s,
then have enough $E$'s so that the resulting $\pi_EC$'s are fat, iterate the
process...  Rather then try to define the precise combinatorics of ``enough",
we take advantage of elementary submodels which provide all the
fat sets we need.

Along with constructing our sequence of functions
$f_{\alpha}$, then, we will construct a sequence of large enough
countable elementary submodels
$\{N_{\alpha}: \alpha < \omega_1\}$ where

\hskip .3 in {\vbox{\hsize 4.5 in

\flushpar\hangindent=.8cm  9.  each $f_{\alpha}, S_{\alpha},
\{N_{\beta}: \beta \leq \alpha\}
\in N_{\alpha+1}$, $\{N_{\beta}: \beta < \alpha\}
\subset N_{\alpha}$, and $\{S_\alpha: \alpha < \omega_1\} \in N_0$.
}}

Further requirements are:

\hskip .3 in {\vbox{\hsize 4.5 in

\flushpar\hangindent=.8cm  1.  if $\alpha \neq \gamma$ then $\{i:
f_{\alpha}(i) =
f_{\gamma}(i)\}$ is finite

\flushpar\hangindent=.8cm  10.  If $C \in N_{\alpha}$ is fat, and $E
= \{n: (n,f_{\alpha}(n)) \in C\}$ then $\pi_EC$ is fat.

\flushpar\hangindent=.8cm  11.  If $C \in N_{\alpha}$ is fat, then
$f_{\alpha} \cap C$ is infinite.

\flushpar\hangindent=.8cm  12.  If $S \in N_{\alpha} \cap [\alpha]^{\omega}$
then $\exists m \ \{\gamma \in S: f_{\gamma} \cap f_{\alpha} \subset m
\times m\}$ is infinite.
}}

Note that property 11 follows from property 10.

As before, property 1 gives us $\{f_{\alpha}: \alpha < \omega_1\}$
almost disjoint.  It remains to show that properties 10 and 11 imply property
8 (which implies property 7, which implies strongly near-Luzin), and property
12 implies there are no uncountable Luzin subfamilies.

\proclaim {Lemma 3.5}  Suppose $\Cal A = \{f_{\alpha}: \alpha < \omega_1\},
\{N_{\alpha}: \alpha < \omega_1\}$ satisfy 1, 9, 10, 11 and 12.  Then

\roster

\item "{(a)}" $\Cal A$ satisfies property 8.

\item "{(b)}"  $\Cal A$ is has no uncountable Luzin subfamily.
\endroster
\endproclaim

\demo {Proof of 3.5 (a)}  We show by induction on $n$ that if
$\beta_0 < ... \beta_n$  and each sup
$S_{\beta_i} <$  $\inf S_{\beta_{i+1}}$ then $\bigcap_{i \leq n}
\bigcup_{\gamma \in S_{\beta_i}}f_{\gamma}$ is fat.
So suppose
$\beta_0, ...
\beta_n$, each sup
$S_{\beta_i} <$  $\inf S_{\beta_{i+1}}$, and $C = \bigcap_{i < n}
\bigcup_{\gamma \in S_{\beta_i}}f_{\gamma}$ is fat.  Let $\alpha <$ inf
$S_{\beta_n}$ with $\beta_{n-1} \in N_{\alpha}$.  Then $C \in
N_{\alpha}$.  Let $\{\gamma_j: j < \omega\} \subset S_{\beta_n}$ with each
$\gamma_j < \gamma_{j+1}$.  Define $\{E_j: j < \omega\}, \{C_j: j < \omega\}$
as follows: $E_{0} = \{k: (k,f_{\gamma_0}(k)) \in C\}$, $C_0 =
\pi_{E_0}C$, $E_{j+1} = \{k: (k,f_{\gamma_{j+1}}(k)) \in C_j\}$, $C_{j+1} =
\pi_{E_{j+1}}C_j$.  By property 9, each $C_j, E_j \in N_{\beta_{j+1}}$.  By
property 10, each
$C_j$ is fat.

But then, by property 11, each $C_j \cap f_{\gamma_j}$ is infinite, so by
$\Cal A$ almost disjoint and $C_j \supset C_{j+1}$, $C \cap \bigcup
\{f_{\gamma_j}: j <
\omega\}$ is fat.

\enddemo

\demo {Proof of 3.5 (b)}  Given an uncountable subfamily
$\Cal B$ of $\Cal A$ and an enumeration  $\Cal B
=\{g_{\alpha}:{\alpha}<{\omega}_1\}$, where
$g_{\alpha}=f_{\phi(\alpha)}$, there are ${\alpha}\le
{\beta}<{\omega}_1$ so
that $S_{\alpha}\subset {\hbox{\rm ran}}\phi\cap {\beta}$ and
$\phi''{\beta}={\hbox{\rm ran}}\phi\cap {\beta}$.
Let $\delta =\phi({\beta})\ge {\beta}$.
Note that ${\alpha} \in N_{\alpha}\subset N_{\beta}\subset
N_{\delta}$, so, by property 12,
$\exists m $ $\{\gamma \in S_{\alpha}: f_{\gamma} \cap f_{\delta} \subset m
\times m\}$ is infinite.  Since $g_{\beta}=f_{\delta}$ and
$\phi^{-1} S_{\alpha}\subset {\beta}$ it follows that
 $\Cal B$ is not Luzin.
\enddemo

Now notice that the construction used in the proof of 2.5 easily adapts to a
construction of a family satisfying properties 1, 9, 11, and 12.  To get
property 10, $r_{k,j}$ is required to satisfy
$$|\pi_{r_{k,j}}C_k| > j,$$ which can be done because $C_k$ is fat.

\bigskip

\noindent {\bf 4.  Trees and anti-Luzin families}
\medskip

The canonical example of an anti-Luzin family is a set of branches of a
countable perfect  tree, i.e. a countable tree such that 
there are two incomparable nodes above every  node.  
What about the reverse?  Must every anti-Luzin
family
look like the branches of a tree?

\definition {Definition 4.1}  An uncountable almost disjoint family $\Cal A$ is
a tree family iff there is a tree ordering $\Cal T = (\bigcup \Cal A, \prec)$
so that for every $a \in \Cal A$ there is a branch $b$ of $\Cal T$ with $a
=^* b$.
\enddefinition

We will show that under CH + $\exists$ a Suslin line there is an anti-Luzin
family which contains no uncountable tree families.

\proclaim {Question 4.2}  Is there (under ZFC alone) an anti-Luzin family
which contains no uncountable tree families?
\endproclaim

While we don't know the answer to question 4.2, we have a related MA +
$\neg$-CH result.

\definition {Definition 4.3}  An almost disjoint family $\Cal A$ is a 
{\it hidden tree family} iff for some infinite $T \subset \bigcup \Cal A$ the set $\{a
\cap T: a \in \Cal A\}$ is a tree family.
\enddefinition

Hidden tree families need not be anti-Luzin.  For example, let
$\Cal A = \{a_{\alpha}: \alpha < \omega_1\}$ be a tree family on the set of
even integers, and let
$\Cal B = \{b_{\alpha}: \alpha < \omega_1\}$ be Luzin on the set of odd
integers.  Then $\{a_{\alpha}
\cup b_{\alpha}: \alpha < \omega_1\}$ is both Luzin and a hidden tree family.

In fact, under MA + $\neg$CH all uncountable almost disjoint families of
size $< 2^{\omega}$ are hidden tree families.

\proclaim {Theorem 4.4}  Assume MA(pre-caliber $\omega_1$).  Then
every uncountable almost disjoint family of
size $< 2^{\omega}$ on $\omega$ is a hidden tree family.
\endproclaim

\demo {Proof of theorem}

We define $p \in {\Bbb P}$:  $p = (T_p,\prec_p,\Cal A_p,h_p)$ where

\roster
\item"{}" $T_p\in [{\omega}]^{<{\omega}}$,  $(T_p,\prec_p)$ is a
finite tree,
\item"{}" $\Cal A_p\subset \Cal A$ is finite,
\item"{}" $h_p:\Cal A_p\to {\omega}$ is a function,
\item"{}" each $(a\setminus h_p(a))\cap T_p$ is linearly ordered by $\prec_p$,
\item"{}" $\forall a\in \Cal A_p$
$\forall n\in a \cap (T_p\setminus h_p(a))$
$\forall k\in T_p\setminus h_p(a)$
if  $k\prec_p n$ then $k\in a$.
\endroster

We define $p \leq q$ iff

\roster
\item"{(a)}" $(T_q,\prec_q)$ is an initial subtree of $(T_p,\prec_p)$,
\item"{(b)}" $\Cal A_q \subset \Cal A_p$,
\item"{(c)}" $h_q\subset h_p$,
\endroster

$\Bbb P$ is easily seen to have pre-caliber $\omega_1$.

\proclaim {Subclaim 4.4.1} For each $a\in \Cal A$, the set
$$
D_a=\{p\in\Bbb P:a\in\Cal A_p \}
$$
is dense in $\Bbb P$.
\endproclaim

\demo {Proof}
Assume that $a\notin \Cal A_p $.
Let $q=(T_p,\prec_p,\Cal A_p \cup \{a\}, h_p\cup(a,n))$,
where $n>\max T_p$. Then  $q\in \Bbb P$ because
$a\cap (T_p \setminus h_p(a))=\emptyset$. The relation $q\le p$ is clear.
\qed\enddemo

\proclaim {Subclaim 4.4.2}  For each $n\in {\omega}$  and $a\in \Cal A$
$$ D_{a,n}=\{p\in \Bbb P:a\in \Cal A_p \ \roman{and} \ a\cap (T_p \setminus n)
\ne\emptyset\}
$$
is dense in $\Bbb P$.
\endproclaim

\demo {Proof}
Let $p\in \Bbb P$. By subclaim 4.4.1 we can assume that $a\in  \Cal A_p $.
Let
$$ k\in a\setminus \bigcup (\Cal A_p\setminus\{a\})\setminus \max \{T_p
\cup n+1\}.$$
Let $q=(T_p\cup \{k\},\prec _q, \Cal A_p, h_p)$, where
$\prec_p\subset \prec_q$ and $\ell\prec_q k$ for each
$\ell\in a\cap (T_p \setminus h_p(a))$. Since $k\notin a'$ for
$a'\in \Cal A_p\setminus \{a\}$  and
$a\cap (T_p\setminus h_p(A)$ is linearly ordered by $\prec_p$
we have that $q\in \Bbb P$ and clearly
$q\le p$.
\qed\enddemo

Let
$$
\Cal D=\{D_{a,n}:a\in \Cal A, n\in {\omega}\}.
$$
By MA($\sigma$-centered)
we have an $\Cal D$-generic filter  $\Cal G$.
Then
$$\Cal T=(T,\prec)
=(\bigcup_{p\in \Cal G}T_p,\bigcup_{p\in \Cal G}\prec_p)
$$
witnesses that $\Cal A$ is a hidden tree family:
taking $h=\bigcup_{p\in \Cal G}h_p$ we have that
$(T\cap a) \setminus h(a)$ is a tail of a branch of $\Cal T$.
\enddemo

In contrast we get
under CH + $\exists$ a Suslin line an anti-Luzin family which has no
uncountable hidden tree families.  In fact we get something stronger.

\definition {Definition 4.5}  An uncountable almost disjoint family $\Cal A$ is
a {\it weak tree family} iff there is a tree ordering $\Cal T = (\bigcup \Cal
A,\prec)$ and a 1-1 function $\phi: \Cal A \rightarrow Br(T)$ (here
$Br(T)$ is the set of branches of $T$) where range
$\phi$ is pairwise disjoint, and each $a \subset^*  \phi(a)$.

$\Cal A$ is a {\it hidden  weak tree family} iff, for some
$T$, $\{a \cap T: a \in \Cal A\}$ is a  weak tree family.
\enddefinition

Weak tree families appeared in \cite{V}  where they are called neat
families. Velickovic proved the following result (lemma 2.3):
{\it Assume $\hbox{MA}_{\aleph_1}$. If $\Cal A\subset [{\omega}]^{\omega}$
is an almost disjoint family
then there is an uncountable family $\Cal B\subset \Cal A$ and a
partition $b=b_0\cup b_1$ for each  $b\in\Cal B$ such that 
$\Cal B_i=\{b_i:b\in\Cal B\}$ is a weak tree family for $i\in 2$}.

\demo{Remark}
One can consider the following weakening of the notion of
weak tree families. 
An uncountable almost disjoint family $\Cal A$ is
a {\it very weak tree family} iff there are a tree ordering 
$\Cal T = (\bigcup \Cal A,\prec)$ and a  function 
$\phi: \Cal A \rightarrow [Br(T)]^{<\omega}$ such that the range
$\phi$ is pairwise disjoint and  each $a \subset^* \bigcup \phi(a)$.
$\Cal A$ is a {\it hidden very weak tree family} iff, for some
$T$, $\{a \cap T: a \in \Cal A\}$ is a very weak tree family.

However,  as it was observed by the referee, 
a hidden very weak tree family can be split into countably many hidden
tree families: for every element $x$ of $\Cal A$ fix a node  of the tree
such that above this node $x$  is covered by a single branch, and spit
$\Cal A$ accordingly.  
\enddemo

\proclaim {Theorem 4.6} Assume CH + $\exists$ a Suslin line.  Then there is an
uncountable anti-Luzin almost disjoint family which contains no uncountable
hidden  weak tree families.
\endproclaim

First, a quick lemma.

\proclaim {Lemma 4.7}  Let $\Cal T^*$ be Aronszajn, $B$ an uncountable set of
branches of $\Cal T^*$ so that no two elements of $B$ have the same order
type.  Then there are incompatible elements
$s, t
\in \Cal T^*$ so that
$\{b \in B: s \in b\}$ and $\{b \in B: t \in b\}$ are uncountable.
\endproclaim

\demo {Proof of lemma 4.7}  By contraposition, suppose $B$ is a set of
branches of $T^*$ with different order types so that if $s,t$ are incompatible
then either
$B(s) =
\{b
\in B: s \in b\}$ is countable or $B(t) =
\{b
\in B: t \in b\}$ is countable.  Then $S = \{s: B(s)$ is
uncountable$\}$ forms a chain, hence is countable.  So there is $\alpha$ with
$T^*(\alpha) \cap S =
\emptyset$, where $T^*(\alpha) =$ the set of elements of $\Cal T^*$ of height
$\alpha$.  But all but countably many elements of $B$ are elements of
$\bigcup_{s \in T^*(\alpha)}B(s)$, so $B$ is countable.
\enddemo

\demo {Proof of theorem 4.6}  Let $\Cal T^*$ be a Suslin tree so that every
element has successors at arbitrarily high levels, and for each
$t
\in \Cal T^*$ construct $c_t \in [\omega]^{\omega}$ so if $s < t$ then $c_s
\supset^* c_t$ and if $s, t$ are not comparable then $c_s \cap c_t =^*
\emptyset$.  Let
$B$ be an uncountable set of branches of $\Cal T^*$ so no two elements of $B$
have the same order type and so that every element of $\bigcup B$ is in
uncountably many branches of $B$.

Let $B = \{b_{\alpha}: \alpha <
\omega_1\}$.

We will define $\Cal A = \{a_{\alpha}: \alpha \in \omega_1\}$ where

 \hskip .3 in {\vbox{\hsize 4.5 in

\flushpar\hangindent=.8cm  13. for all
$s \in b_{\alpha}, c_s \supset^* a_{\alpha}$.
}}

$\Cal A$ will clearly be almost disjoint.

Let $\{\Cal T_{\alpha} = (T_{\alpha}, \prec_{\alpha}): \alpha < \omega_1\}$
enumerate all  perfect trees whose underlying set is some infinite subset of
$\omega$.

We further require

 \hskip .3 in {\vbox{\hsize 4.5 in

\flushpar\hangindent=.8cm  14.  for all $\beta < \alpha$ either for some $s \in
b_{\alpha} \ c_s \cap T_{\beta}$ is contained, mod finite, in
a branch of
$\Cal T_{\beta}$, or $a_{\alpha}\cap T_{\beta}$ is not a subset, mod
finite,  of a branch of $\Cal T_{\beta}$. }}

The family $\Cal A$ is constructed recursively in ${\omega}_1$ many steps. In
the ${\alpha}^{\hbox{th}}$ step we apply lemma 4.8 below to get $a_{\alpha}$.

\proclaim {Lemma 4.8}  Suppose $\{a_{\beta}:
\beta < \alpha\}$ satisfies property 13.  Then there is a set $a$ satisfying
properties 13 and 14.
\endproclaim

\demo {Proof}  Rather than describe the proof as an induction, we will
(equivalently) use the Rasiowa-Sikorski lemma
(see \cite{K,Theorem 2.21}), defining a countable set of
forcing conditions and countably many dense sets so that any generic filter
meeting the dense sets gives rise to the desired object.

The partial order is as follows:  ${\Bbb P}$ consists of all pairs $p =
(a_p,b_p)$ where $a_p$ is a finite subset of
$\omega$ and $b_p$ is a
finite subset of $b_{\alpha}$.  The order is as follows:  $p \leq q$ iff $a_p
\supset a_q$,
 $b_p \supset b_q$, and $a_p
\setminus a_q \subset \bigcap_{t \in b_q}c_t$.

Clearly $\{p: |a_p| > n\}$ and $\{p: s \in b_p\}$ are dense for each $n <
\omega, s \in b_{\alpha}$, so if
$\Cal G$ is a filter meeting each of these dense sets then $\bigcup_{p \in \Cal
G}a_p$ satisfies property 13.

Towards property 14, fix $\Cal T = \Cal T_{\beta}, T = T_{\beta}$.  We may
assume that for every
$t
\in b_{\alpha} \ c_t \cap T$ is not a subset, mod finite, of a branch
of $\Cal T$.  For each $n < \omega$ define  
$$D(\Cal T,n) = \{p: a_p\setminus n \
\hbox{contains  $\Cal T$-incomparable elements} \}.$$  
We show that $D(\Cal T,n)$ is dense for each $n$.

Fix
$q \notin D(\Cal T,n)$.  Let $c = \bigcap_{t \in b_q}c_t$.  Since $c \cap T
\setminus n$ is not a subset, mod finite, of a
branch of $\Cal T$, there are two  $\prec_{\beta}$-incompatible
elements, ${\xi}$ and ${\eta}$, of
$T \cap (c \setminus n)$.  Set $a_p = a_q\cup\{{\xi},{\eta}\}$, 
$b_p = b_q$.  Then $p \in D(\Cal T,n)$ and $p \leq q$.

If, for all $n$, $\Cal G$ meets $D(\Cal T,n)$, $\bigcup_{p \in \Cal G}a_p \cap
T$ will not be a subset, mod finite,  of any  branch of $\Cal T$.
Lemma 4.8 is proved.

The
following two lemmas, once proved,  will
complete the proof of Theorem 4.6.

\proclaim {Lemma 4.9}  If property 14 holds, $\Cal A$ has no uncountable
hidden weak tree families.
\endproclaim

\demo  {Proof} This is where we use that $\Cal T^*$ is Suslin.

So suppose $\Cal B$ is an uncountable subset of $\Cal A$, $T \subset
\omega$ is infinite, and $\Cal C = \{a
\cap T: a
\in \Cal B\}$ is a collection of infinite sets.  We show that $\Cal C$ is not a
very weak tree family.

Let $\Cal T = (T, \prec)$ and $\phi: \Cal C \rightarrow Br(\Cal
T)$.  We show that $\phi$ does not satisfy the properties of
definition 4.5.

For some $\beta, \Cal T = \Cal T_{\beta}$.

Let $S = \{s \in \Cal T^*: c_s \cap T$ is a subset, mod finite,  of a 
branch of $\Cal T\}$.  If $S = \emptyset$, then by property 14 for all
but countably many $a \in \Cal B$, $a \cap T \not\subset^* \phi(a
\cap T)$.

Suppose $s \in S$, $c_s \cap T \subset^* d$, where  $d$
is a branch of $\Cal T$.  Then if $s \in b_{\alpha}$, $a_{\alpha} \cap T
\subset^* d$.  Let $S'$ be the set of $\Cal T^*$-minimal elements of
$S$.  By $\Cal T^*$ Suslin, $S'$ is countable.  So there is $s \in S'$ with
$\{b_{\alpha}: a_{\alpha} \in \Cal B$ and $s \in b_{\alpha}\}$ uncountable.
But then either there are uncountably many $\alpha$ with 
$a_{\alpha} \not\subset
\phi(a_{\alpha} \cap T)$, or $\phi$ is not 1-1.
\enddemo

\proclaim {Lemma 4.10}  If property 13 holds, $\Cal A$ is anti-Luzin.
\endproclaim

\demo {Proof}  Suppose $\Cal B$ is an uncountable subset of $\Cal A$.  By
lemma 4.7 there are incompatible
$s, t \in \Cal T^*$ with $\Cal C = \{b_{\alpha} \in \Cal B: s \in b_{\alpha}\}$
and
$\Cal D =
\{b_{\alpha} \in \Cal B: t \in b_{\alpha}\}$ uncountable.  Without loss
of generality, we may assume that for some $n \geq$ sup $c_s \cap c_t$ if
$b_{\alpha}
\in
\Cal C$ then $a_{\alpha} \setminus c_{s} \subset n$ and if $b_{\alpha} \in \Cal
D$ then
$a_{\alpha}
\setminus c_{t} \subset n$.  But
then
$\bigcup \Cal C \cap \bigcup \Cal D \subset n$, as desired.
\enddemo
\enddemo
\enddemo
\bigskip

\noindent {\bf 5.  Between near-Luzin and strongly near-Luzin}
\medskip

\definition {Definition 5.1}  An uncountable almost disjoint family $\Cal A$ is
$k$-near-Luzin iff for every $\Cal C_0,... \Cal C_{k-1} \in
[\Cal A]^{\omega_1}$, $\bigcap_{i<k}\bigcup \Cal C_i$ is infinite.
\enddefinition

The purpose of this section is to show that these notions are
(consistently) distinct.

Clearly near-Luzin is $2$-near-Luzin and so every Luzin family is
$2$-near-Luzin, but not necessarily contains $3$-near-Luzin subfamily, as
we will see it in theorem 5.5.

\proclaim {Theorem 5.2}  The following is consistent:  $\forall k \in
[2,\omega)$ there is an uncountable almost disjoint family $\Cal A_k$ so $\Cal
A_k$ is
$k$-near-Luzin, contains no uncountable Luzin subfamilies, and contains no
uncountable
$k+1$-near-Luzin subfamilies.
\endproclaim

The proof proceeds by showing that for every $k$ there is a partial order
${\Bbb P}_k$ with precaliber $\omega_1$ forcing $\Cal A_k$ to exist, and
iterating with precaliber $\omega_1$.  It is easy to see that both $k$-near
Luzin and ``no Luzin subfamilies" are preserved by precaliber $\omega_1$
forcing.  The way we ensure no $k+1$-near-Luzin subfamilies will also be
preserved by precaliber $\omega_1$ forcing.

Unlike our earlier constructions, each $|\Cal A_k| = 2^{\omega}$.  Each
$\Cal A_k$ is again a family of functions, but instead of functions on
$\omega$ the domains come from a $k$-linked not ($k+1$)-linked family
$\Cal E_k$ with
special properties.  This family was first constructed by Hajnal; the
construction appeared in \cite{JS}.

\proclaim {Lemma 5.3}  For all $k < \omega$ there is a family 
$\Cal E_k\subset [{\omega}]^{\omega}$ with
$|\Cal E_k| = 2^{\omega}$ so that

\hskip .3 in {\vbox{\hsize 4.5 in

\flushpar\hangindent=.8cm 15. if $e_0,... e_{k-1} \in \Cal E_k$ then
$\bigcap_{i<k}e_i$ is infinite

\flushpar\hangindent=.8cm 16.  if $e_0,... e_k$ are distinct elements
of $\Cal E_k$ then $|\bigcap_{i\leq k}e_i| < \omega$

\flushpar\hangindent=.8cm 17.  if $X \in [\Cal E_k]^{\omega_1}$ there
are $X_0,... X_k \in [X]^{\omega_1}$ with $$|\bigcap_{i \leq
k}\bigcup X_i| < \omega.$$
}}
\endproclaim

\demo {Proof}  
Let $S_k=\{\left[{}^n2\right]^k:n<{\omega}\}$. 
We will construct $\Cal E_k$ as a subset of $[S_k]^{\omega}$ but 
since $|S_k|={\omega}$, this construction proves the lemma.
For $f: \omega \rightarrow 2$ let $e_f =
\bigcup_{n<\omega}\{s
\in [2^n]^k: f{\restr}n \in s\}$.  Let $\Cal E_k = \{e_f: f \in 2^{\omega}\}$.

We show that property 15 holds:  Given $f_0,... f_{k-1}$ distinct, pick
$m$ so
$\{f_i{\restr}m: i < k\}$ are distinct.  But then for each $n \geq m$
$$\{f_i{\restr}n: i < k\} \in \bigcap_{i<k}e_{f_i}.$$

We show that property 16 holds:  Given $f_0,... f_{k}$ distinct, there
is
$m$ so
$\{f_i{\restr}m: i \leq k\}$ are distinct.  But then $\bigcap_{i \leq
k}e_{f_i} \subset \bigcup_{j<m}[2^j]^{k+1}$.

We show that property 17 holds:  Fix $X \in [\Cal E]^{\omega_1}$.  Let
$I = \{f: e_f \in X\}$ and let
$g_0,... g_k$ be distinct complete accumulation points of $I$ in the
usual topology on $2^{\omega}$.  Fix $n$ so $g_0{\restr}n,... g_k{\restr}n$ are distinct.
Define $X_i = \{e_f \in \Cal E: f{\restr}n = g_i{\restr}n\}$.  Then $$\bigcap_{i
\leq k}\bigcup X_i \subset
\bigcap_{j<n}[2^j]^k$$ which completes the proof.

\enddemo

Note that by construction properties 15 and 16 are absolute in the
following sense:  let $\Cal E = (\Cal E_k)^M$, and let $M \subset N$,
where $M, N$ are models of enough set theory.  Then 15 and 16 hold for $\Cal E$
in
$N$.

The next lemma says that property 17 is preserved in some models. In this and
succeeding proofs we will refer to the following easy fact about ccc forcing:

\proclaim {Fact 5.4}  If ${\Bbb P}$ is a ccc partial order, $P$ is an
uncountable  subset of ${\Bbb P}$ and $\dot G$ names the generic
filter, then there is $p \in {\Bbb P} \ p \Vdash P\cap \dot G$ is
uncountable.
\endproclaim

\proclaim {Lemma 5.5}  Suppose $\Cal E$ has the property that if $X \in
[\Cal E]^{\omega_1}$ there are $X_0,... X_k \in [X]^{\omega_1}$ with
$$|\bigcap_{i \leq k}\bigcup X_i| < \omega.$$  Then $\Cal E$ will still
have this property in a forcing extension by a precaliber $\omega_1$
partial order.
\endproclaim

\demo {Proof}  Let ${\Bbb P}$ have precaliber $\omega_1$, and suppose
$\Vdash_{\Bbb P} \dot X = \{\dot x_{\alpha}: \alpha < \omega_1\} \in [\Cal
E]^{\omega_1}$.  Fix $p \in {\Bbb P}$.  For each $\alpha$ pick $p_{\alpha}
\leq p$ so for some
$e_{\alpha} \in \Cal E \ p_{\alpha} \Vdash \dot x_{\alpha} =
e_{\alpha}$.  Then there is an uncountable centered family
$\{p_{\alpha}:
\alpha \in I\}$.  Let $Y = \{e_\alpha: \alpha \in I\}$.  By hypothesis
there are $I_0,... I_k \in [I]^{\omega_1}$ and $m < \omega$ with
$$\bigcap_{i\leq k}\bigcup\{e_{\alpha}: \alpha \in I_i\} \subset m.$$

Let $\dot G$ be the generic filter, and define $$\dot J_i = \{\alpha
\in I_i: p_\alpha \in \dot G\}.$$  List each $I_i$ as $\{\alpha^i_{\gamma}:
\gamma < \omega_1\}$.  Let $q_{\gamma} \leq
p_{\alpha^i_{\gamma}}$ for all $i \leq k$.  Let $Q = \{q_{\gamma}: \gamma <
\omega_1\}$.  By fact 5.4 there is $p \in {\Bbb P}$ with $p \Vdash |Q \cap
\dot G| = \omega_1$.  So $\forall i \leq k \ p \Vdash \dot J_i$ is
uncountable, which by a density argument completes the proof.

\enddemo

\proclaim {Lemma 5.6}  Suppose $\Cal E \subset \wp(\omega)$ satisfies the
following:

\hskip .3 in {\vbox{\hsize 4.5 in

\flushpar\hangindent=.8cm  if $X \in [\Cal E]^{\omega_1}$ there
are $X_0,... X_k \in [X]^{\omega_1}$ with $$|\bigcap_{i \leq
k}\bigcup X_i| < \omega.$$
}}

If $\Cal A = \{f_e: e \in \Cal E\}$ where each $f_e: e \rightarrow \omega$
then $\Cal A$ has no uncountable ($k+1$)-near Luzin subfamilies.

\endproclaim

\demo {Proof}  If for some finite m$$\bigcap_{i \leq
k}\bigcup X_i \subset m$$ then $$\bigcap_{i \leq k}\bigcup_{e \in
X_i}f_e \subset m \times \omega.$$  Let $Y_i \in [X_i]^{\omega_1}$ so $\exists
\sigma_i \ \forall e
\in Y_i \ f_e{\restr}m = \sigma_i$.  Then $$\bigcap_{i \leq k}\bigcup_{e \in
Y_i}f_e \subset \bigcap_i \sigma_i.$$
This completes the proof of lemma 5.6.
\enddemo

\demo{Proof of theorem 5.2}
Let $\Cal E_k\subset [{\omega}]^{\omega}$ be a family satisfying
15--17
from
lemma 5.3. We will have $\Cal A_k=\{f_e:e\in \Cal E_k\}$, where 
$f_e:e\to {\omega}$ is a function for $e\in \Cal E_k$.
By lemma 5.6, this assumption guarantee that $\Cal A_k$
has no $(k+1)$-near-Luzin subfamily.  

Define ${\Bbb P}_k$, a pre-caliber $\omega_1$ forcing which adds
generic almost disjoint functions $f_e:e\to {\omega}$ for $e\in \Cal E_k$
as follows:
\roster

\item"{}" elements of ${\Bbb
P}_k$ have the form $p = \{\sigma_{p,e}: e \in E_p\}$,

\item"{}" $E_p$ is a finite
subset of
$\Cal E_k$,

\item"{}" each $\sigma_{p,e}$ is a finite function from $e$ to $\omega$.
\endroster

The
order is:
$p
\leq q$ iff

\roster

\item"{}" $E_p \supset E_q$

\item"{}" for $e \in E_q \ \sigma_{p,e} \supset
\sigma_{q,e}$

\item"{}" for $e \neq e' \in E_q \ \sigma_{p,e} \cap \sigma_{p,e'}
= \sigma_{q,e} \cap \sigma_{q,e'}$.
\endroster

${\Bbb P}_k$ is easily seen to have pre-caliber
$\omega_1$.  We define $p
\Vdash \dot f_e(i) = j$ iff $e \in E_p$ and $\sigma_{p,e}(i) = j$.  By a
standard genericity argument, $\Vdash_{{\Bbb P}_k}$ [ $\dom \dot f_e = e$ and
if $e
\neq e'$ then $|\dot f_e \cap \dot f_{e'}| < \omega $].

\proclaim {Lemma 5.7}  $\Vdash_{{\Bbb P}_k} \dot A_k$ is $k$-near-Luzin.
\endproclaim

\demo {Proof}  Working in $V^{{\Bbb P}_k}$, suppose for each $i < k$ we have
$\dot X_i$ an uncountable subset of $\Cal E_k$.  We want to show that
$\bigcap_{i<k}\bigcup \{\dot f_e: e \in \dot X_i\}$ is infinite.

We may assume the $\dot X_i$'s are disjoint, and each $\dot X_i = \{\dot
e_{\alpha,i}:
\alpha <
\omega_1\}$ in a 1-1 enumeration.  Fix $p \in {\Bbb P}_k$.  For each $\alpha$
there is
$p_{\alpha}
\leq p$ and for each
$\alpha,i$ there is
$d_{\alpha,i}$ so $$p_{\alpha} \Vdash \forall i < k \ \dot e_{\alpha,i} =
d_{\alpha,i}$$ and for each $\alpha$ the $d_{\alpha,i}$'s are distinct.

We may assume the $p_{\alpha}$'s are centered.

Since the $p_{\alpha}$'s are centered, and the enumeration is 1-1,
$$\forall i \ \roman{if} \ \alpha \neq \beta \ \roman{then}
\ d_{\alpha,i} \neq d_{\beta,i}.$$

Pick distinct $\alpha_0,... \alpha_{k-1}$ so that for $i \neq j \
d_{\alpha_i,i} \neq d_{\alpha_j,j}$ and (by a $\Delta$-system
argument) $d_{\alpha_i,i}
\notin E_{p_{\alpha_j}}$.  By property 15
$$|\bigcap_{i<k} d_{\alpha_i,i}| = \omega.$$  Let $m \in\bigcap_{i<k}
d_{\alpha_i,i}, m >$ sup  $\dom \sigma_{p_{\alpha_i},e}$
for all $i$, all $e \in E_{p_{\alpha_i}}$.  There is $q < p_{\alpha_i}$ for
all $i$ with $\sigma_{q,d_{\alpha_i,i}}(m) = 0$.  So
$q \Vdash \exists m > n (m,0) \in
\bigcap_{i<k}\bigcup \{\dot f_e: e \in \dot X_i\}$. A
density argument completes the proof.
\enddemo

\proclaim {Lemma 5.8} $\Vdash_{{\Bbb P}_k} \dot A_k$ has no uncountable Luzin
subfamilies.
\endproclaim

\demo {Proof}  Suppose $\{\dot a_{\dot e_{\alpha}}: \alpha < \omega_1\} \subset
\Cal A_k$ and $p \Vdash$ the enumeration $\{\dot a_{\dot e_{\alpha}}: \alpha <
\omega_1\}$ witnesses that the family is Luzin.  Choose $p_{\alpha} \leq p,
d_{\alpha}
\in
\Cal E_k$ with $p_{\alpha} \Vdash \dot e_{\alpha} = d_{\alpha} \in
E_{p_{\alpha}}$.  We may assume

\hskip .3 in {\vbox{\hsize 4.5 in

\flushpar\hangindent=.8cm  18.  the
$p_{\alpha}$'s are centered

\flushpar\hangindent=.8cm  19.  $\{E_{p_{\alpha}}: \alpha < \omega\}$ is a
$\Delta$-system with root $E$

\flushpar\hangindent=.8cm  20.  $\exists n \forall \alpha \ E_{p_{\alpha}} =
\{e_{\alpha,i}: i < n\}$

\flushpar\hangindent=.8cm  21.  $\forall i \exists \sigma_i \forall \alpha
\sigma_{p_{\alpha},e_{\alpha,i}} = \sigma_i$
}}

By necessity

\hskip .3 in {\vbox{\hsize 4.5 in

\flushpar\hangindent=.8cm  22.  $\exists m \forall i \ \sigma_i \subset m
\times m$.
}}

There is $q \leq p_{\omega}$ and $k$ such that $q \Vdash \forall i > k \ \dot
a_{\dot e_i} \cap \dot a_{\dot e_{\omega}} \not\subset m \times m$.

By property 19 $\exists j \geq k$ with $E_q \cap E_{p_j} = E$.

We define $r \leq q$:

\hskip .3 in {\vbox{\hsize 4.5 in

\flushpar\hangindent=.8cm  $E_r = E_q \cup E_{r_j}$

\flushpar\hangindent=.8cm  for $e \in E_q$, $\sigma_{r,e} = \sigma_{q,e}$

\flushpar\hangindent=.8cm  for $e \in E_{r_j} \setminus E$, $\sigma_{r,e} =
\sigma_{r_j,e}$.
}}

Then $r \Vdash \dot a_{\dot e_j} \cap \dot a_{\dot e_{\omega}} = \dot a_{d_j}
\cap \dot a_{d_{\omega}} \subset m \times m$, a contradiction.

\enddemo

Theorem 5.2 is proved.
\enddemo
Finally, we note that Luzin does not imply $3$-near Luzin

\proclaim {Theorem 5.9}  There is a Luzin almost disjoint family with no
uncountable
$3$-near Luzin subfamily.
\endproclaim

\demo {Proof}  Let $\Cal E = \Cal E_2$ be as in lemma 5.3.  As in theorem 5.2,
we construct $\Cal A =
\{f_e: e
\in \Cal E\}$ where each  $\dom f_e = e$, so $\Cal A$ has no uncountable
$3$-near Luzin subfamily.  Here is how we get Luzin.

Let $\Cal E = \{e_{\alpha}: \alpha < \omega_1\}$, $f_{\alpha} =
f_{e_{\alpha}}$.

Our induction hypothesis
at stage $\alpha$ is that for all $\beta < \alpha$ and all $n < \omega
\ \{\gamma < \beta: f_{e_{\gamma}} \cap f_{e_{\beta}} \subset n \times
\omega\}$ is finite.  This will certainly give us Luzin.

At stage $\alpha$ fix a 1-1 enumeration $\{\beta_n: n < \omega\}$ of
$\alpha$.  In the $n^{\roman{th}}$ step of the construction of $f_{\alpha}$
we ensure that  $\dom f_{\alpha} \cap n = e_{\alpha} \cap n$ and $f_{\alpha}
\cap f_{\beta_n} \not\subset n \times \omega$, without increasing
$f_{\alpha} \cap f_{\beta_m}$ for $m < n$.  Since $e_{\alpha} \cap
e_{\beta_n} \setminus \bigcup\{e_{\beta_m}: m < n\}$ is finite, this can be
done, and the construction is complete.
\enddemo


\bigskip
\bigskip
\bigskip
\bigskip
\bigskip
\bigskip
\bigskip
\bigskip

{\leftskip=8cm
Nov1997
\medskip

}

\enddocument